\documentclass[12pt]{article}
\usepackage{amssymb, amsmath, amsthm}
\usepackage{hyperref}

\newcommand{\ex}{{\rm ex}}
\newcommand{\Prob}{{\rm Prob}}
\newcommand{\Perm}{{\rm Perm}}
\newcommand{\Card}{{\rm Card}}
\newcommand{\A}{{\mathfrak A}}
\newcommand{\Nat}{{\mathbb N}}
\newtheorem{theorem}{Theorem}

\newtheorem{corollary}[theorem]{Corollary}
\newtheorem{lemma}[theorem]{Lemma}

\theoremstyle{definition}
\newtheorem{definition}[theorem]{Definition}
\newtheorem{remark}[theorem]{Remark}

\begin{document}

\title{The  boundary of the  Eulerian number triangle}
\author{Alexander Gnedin\footnote{Utrecht University,
Mathematisch Instituut, PO Box 80010, 3508 TA Utrecht, The Netherlands;
gnedin@math.uu.nl\,.}\,\,\,\, and \setcounter{footnote}{6} Grigori
Olshanski\footnote{Institute for Information Transmission Problems, Bolshoy
Karetny 19,  Moscow 127994, Russia; olsh@online.ru\,. Research of G.O. was
supported by the CRDF grant RUM1--2622--ST--04.}}

\date{}

\maketitle

\begin{abstract}
\noindent The Eulerian triangle is a classical array of combinatorial numbers
defined by a linear recursion. The associated boundary problem asks one to find
all extreme nonnegative solutions to a dual recursion. Exploiting connections
with random permutations and Markov chains we show that the boundary is
discrete and explicitly identify its elements.
\end{abstract}

\section{Introduction and main results}
The {\it Eulerian triangle\/} (see, e.g., \cite[Section 6.2]{GKP}) is the
infinite array of {\it Eulerian numbers} \footnote{ Other commonly used
notations  are $A_{nk}$ and $E_{nk}$. Some authors follow historical definition
of Eulerian numbers as coefficients of Eulerian polynomials, hence index the
numbers by $k$ ranging from $1$ to $n$.}
$$
\left\langle {n \atop k} \right\rangle\, \qquad (0\le k\le n-1\,, ~n=1,2,\ldots),
$$
considered together with the defining them recursion
\begin{equation}\label{dim}
\left\langle {n \atop k} \right\rangle = (k+1) \left\langle {n-1 \atop k}
\right\rangle + (n-k)\left\langle {n-1 \atop k-1} \right\rangle  \,\qquad
\end{equation}
and the boundary conditions
\begin{equation}\label{bc}
\left\langle{1 \atop 0} \right\rangle=1\,,~~~
{\rm and~~}\left\langle {n \atop k}
\right\rangle = 0 ~~~ \text{for\, $n<0$\, or\, $k>n-1$}.
\end{equation}
The first six rows of the triangle ($n=1,\dots,6$) are
\begin{eqnarray*}
\begin{array}{ccccccccccc}
                            &&&&&1&&&&&\\
                            &&&&1&&1&&&&\\
                            &&&1&&4&&1&&&\\
                            &&  1 && 11&& 11&& 1&&\\
                             & 1 & & 26 & & 66 & & 26 & & 1&\\
                             1&& 57&& 302&& 302&&57&& 1\\
\end{array}
\end{eqnarray*}

\par
We are interested in  {\it nonnegative\/} solutions
$$
V=(V_{nk}:\quad 0\le k\leq n-1\,, \quad n=1,2,\dots)
$$
to the {\it dual\/} or {\it backward} recursion
\begin{equation}\label{dual}
V_{nk}= (k+1)V_{n+1,k}+ (n-k) V_{n+1,k+1},
\end{equation}
subject to the only normalization condition $V_{10}=1$. In contrast to
(\ref{dim}), the dual recursion has multiple solutions which comprise a convex
set $\cal V$. We denote by $\ex({\cal V})$ the set of extreme solutions and
call it the {\it  boundary}.

\par
One principal result of this paper gives a natural parametrization of the
boundary:

\vskip0.5cm
\begin{theorem}\label{thm1}
Every extreme solution $W=(W_{nk})\in {\rm ex}({\cal V})$ is uniquely determined
by the parameter $\theta:=W_{20}$, which assumes
values in the following subset of the unit interval
$$
\Theta=\left\{
\frac12\cdot\frac{\varkappa}{\varkappa+1}\,,\quad\varkappa=0,1,\ldots\right\}
\cup \left\{\frac12\right\}\cup
 \left\{ \frac12\cdot\frac{\varkappa+2}{\varkappa+1}\,,\quad\varkappa=0,1,\ldots\right\}.
$$
The correspondence between $\ex(\cal{V})$ and $\Theta$ is a homeomorphism.
\end{theorem}

\noindent Thus, the parameter set $\Theta$ is composed of two sequences and
their sole accumulation point $\frac12$. An obvious symmetry of $\Theta$ about
$\frac12$ corresponds to  the symmetry of (\ref{dim}) under the substitution
$(n,k)\to (n,n-1-k)$.

\par
Theorem \ref{thm1} distinguishes the Eulerian triangle from other classical
number triangles, whose boundaries can be also identified with subsets of the
unit interval. It is well known that the boundary of the Pascal triangle is
$[0,1]$; a result equivalent to both de Finetti's theorem on exchangeable
trials and Hausdorff's characterization of moment sequences \cite{Aldous,
Blackwell, Kerov.young, Sawyer}. The boundary of the $q$-Pascal triangle of
$q$-binomial coefficients is $\{0\}\cup\{\ldots, q^{-3}, q^{-2}, q^{-1}\}$ for
$q>1$, see \cite{Kerov.book, Olshanski1}. For a parametric family of the
generalized Stirling  triangles  considered in \cite{Gibbs} the boundary was
shown discrete (as for the $q$-Pascal triangle) for some values of the
parameter and coinciding with $[0,1]$ for other. For instance, the boundary is
discrete for the triangle of Stirling numbers of the second kind, and it is
continuous for the triangle of signless Stirling numbers of the first kind.
Some sporadic earlier results on a wider family of generalized Stirling
triangles are found in  \cite{Kerov.af}.

\par
The extreme solution corresponding to a particular value of the parameter
$\theta\in\Theta$ will be denoted $W(\theta)=(W_{nk}(\theta))$. Our second
principal result gives explicit formulas for  these solutions:
\vskip0.5cm
\begin{theorem}\label{thm2}
For
$\theta=\frac{\varkappa+2}{2(\varkappa+1)}>\frac12$ with $\varkappa=0,1,2,\dots$
we have
\begin{equation}\label{Wnk1}
W_{nk}(\theta)= {n+\varkappa-k\choose n}{\bigg/}(\varkappa+1)^n \quad (
n=1,2,\dots;\, 0\le k\le n-1),
\end{equation}
whereas for  $\theta=\frac{\varkappa}{2(\varkappa+1)}<\frac12$
with $\varkappa=0,1,2,\dots$ we have
\begin{equation}\label{Wnk2}
W_{nk}(\theta)= {\varkappa+k+1\choose n}{\bigg/}(\varkappa+1)^n \quad (
n=1,2,\dots;\, 0\le k\le n-1).
\end{equation}
Finally, the solution corresponding to $\theta=\frac12$ is
\begin{equation}\label{Wnk3}
W_{nk}(\tfrac12)= \frac1{n!} \quad( n=1,2,\dots;\, 0\le k\le n-1).
\end{equation}
\end{theorem}

\vskip0.5cm
\noindent
It is seen that for $\theta>\frac12$,
$W_{nk}(\theta)$ is $0$ for $k>\varkappa$, and for
$\theta<\frac12$ it is $0$ for
$k< n-1-\varkappa$. Notice also the
symmetry
$$
W_{nk}(\theta)=W_{n,n-1-k}(1-\theta).
$$
All three formulas of Theorem \ref{thm2} can also be written in a unified way
\begin{equation}\label{GrishinaF}
W_{nk}(\theta)=\frac1{n!}\,\prod_{i=-k}^{-k+n-1}(1+\theta'i),
\end{equation}
where $\theta':=2\theta-1$ ranges over the set
$$
\Theta'=\{-1,\,-\tfrac12,\,-\tfrac13,\,-\tfrac14,\,\dots,\,
0,\,\dots,\,\tfrac14,\,\tfrac13,\,\tfrac12,1\}.
$$

\par
Theorems \ref{thm1} and \ref{thm2} are our main results. They imply a simple
description of the whole set of nonnegative normalized solutions to (\ref{dual}).
\vskip0.5cm
\begin{corollary}\label{uniqueness}
Each
solution $V\in\cal V$
can be uniquely represented as a convex combination
$$
V=\sum_{\theta\in\Theta}p(\theta)  \,W(\theta)\,,
$$
with $p$ a probability distribution on $\Theta$.
\end{corollary}

\par
The rest of the paper is organized as follows. In Section \ref{M} we encode the
Eulerian triangle into a graded graph $\cal E$, which allows to embed our
boundary problem in the general framework developed in \cite{Kerov.af,
Kerov.book, KOO, VK-FA}. In particular, we recall a construction of the Martin
boundary, which is our major technical tool. In Section \ref{D} we explain
connections to  random permutations and random total orders on
$\Nat=\{1,2,\dots\}$ with a sufficiency property. In Section \ref{B} the
extreme solutions are related to random sorting algorithms. The main results
are proved in  Sections \ref{P} and \ref{EOP}. Some connections and extensions
are discussed in the last section.

\par
Boundary problems of combinatorial nature
arise within a variety of mathematical disciplines like probability, numerical
analysis, algebra and representation theory. Our own view on the subject has
two major sources: the asymptotic representation theory of symmetric groups
initiated by Kerov and Vershik \cite{VK-FA} and Kingman's theory of
exchangeable partitions summarized in  lectures by Pitman \cite{CSP}.

\section{The graph $\cal E$ and the Martin boundary}
\label{M}

The graph $\cal E$ is an infinite directed graded graph with the vertex set
$$\{(n,k); ~n=1,2,\dots,~ 0\leq k\leq n-1\}$$
and multiple edges. The  $n$th level of $\cal E$, denoted ${\cal E}_n$, is the
set of $n$ vertices $(n,0),\ldots,(n,n-1)$. The level ${\cal E}_1$ has a single
vertex $(1,0)$, which  is the root of $\cal E$. For generic vertex $(n,k)$ the
outgoing edges link $(n,k)$ to two vertices $(n+1,k)$ and $(n+1,k+1)$ by $k+1$
and $n-k$ directed edges, respectively. We say that two vertices connected by
an edge are {\it adjacent} with each other. A path in $\cal E$ (finite or
infinite) starting at $(n,k_n)$  is a sequence of edges linking adjacent
vertices $(n,k_n), (n+1,k_{n+1}),\ldots$ on consecutive levels. A path starting
at the root $(1,0)$ is called {\it standard}. The  edge multiplicities in $\cal
E$ are selected to  match with the  coefficients in the recursion \eqref{dim},
from which it is clear that the Eulerian number $\left\langle {n \atop k}
\right\rangle$ is the {\it dimension} \footnote{The dimensions alone do not
determine $\cal E$ unambiguously, since the Eulerian numbers  satisfy  many
other recursions different from \eqref{dim}. By our definition of the Eulerian
triangle we were careful to say that  the array was considered together with
the recursion \eqref{dim}.} of the vertex $(n,k)$, meaning the number of
standard paths in $\cal E$ which terminate at $(n,k)$.

\par
The interpretation of dual recursion (\ref{dual}) requires some concepts of
probability theory. Let us consider $\cal E$ as a state space of some Markov
process whose time parameter $n$ runs in the {\it reverse\/} direction
$\dots\to n\to n-1\to\dots\to1$, a possible state at time $n$ is a vertex in
${\cal E}_n$, and the transition probabilities are given by
\begin{gather}
\Prob\{(n,k)\to (n-1,k)\}=(k+1)\frac{\left\langle {n-1 \atop k}
\right\rangle}{\left\langle {n \atop k} \right\rangle}\label{prob1}\\
\Prob\{(n,k)\to(n-1,k-1)\}=(n-k)\frac{\left\langle {n-1 \atop k-1}
\right\rangle}{\left\langle {n \atop k} \right\rangle}\,.\label{prob2}
\end{gather}
The basic relations \eqref{dim} and \eqref{bc} translate as the rule
of total probability

\begin{equation}\label{sumprob}
\Prob\{(n,k)\to (n-1,k)\}+\Prob\{(n,k)\to(n-1,k-1)\}=1,
\end{equation}
and imply that at consecutive times the process must reside in adjacent
vertices of $\cal E$.

\par
Now let $V=(V_{nk})\in {\cal V}$.
Setting
\begin{equation}\label{tilde}
\widetilde V_{nk}=\left\langle {n \atop k} \right\rangle V_{nk},
\end{equation}
the recursion \eqref{dual} can be rewritten as
\begin{equation}\label{dual2}
\begin{split}
\widetilde V_{nk}&=\Prob\{(n+1,k)\to(n,k)\} \tilde V_{n+1,k}
\\&
+\Prob\{(n+1,k+1)\to(n,k)\} \widetilde V_{n+1,k+1}
\end{split}
\end{equation}

\begin{lemma}\label{norm}
We have
\begin{equation}\label{normalisation}
\sum_{k=0}^{n-1}\widetilde V_{nk}=1, \qquad n=1,2,\dots.
\end{equation}
\end{lemma}
\begin{proof}
Indeed, from \eqref{sumprob} and \eqref{dual2}, the quantity
$\sum_{k=0}^{n-1}\widetilde V_{nk}$ does not depend on $n$. Since it equals $1$
for $n=1$ due to the normalization condition, the same holds for all $n$.
\end{proof}
\noindent Thus, the vector $(\widetilde V_{n0},\dots,\widetilde V_{n,n-1})$ is
a probability distribution on ${\cal E}_n$ for each $n$, and this family of
distributions is consistent with the transition probabilities \eqref{prob1} and
\eqref{prob2}. It follows that $V$ determines the law of a Markov chain by the
virtue of \eqref{tilde}. The boundary problem acquires therefore the following
meaning:
\begin{itemize}
\item
describe all possible probability laws for a Markov chain on $\cal E$,
whose transition probabilities are given by \eqref{prob1} and \eqref{prob2}.
\end{itemize}

\par
If we required (\ref{dual}) to only hold for $n$ restricted to some finite
range $1\le n\le N$, the analogous boundary problem were rather simple. For,
each truncated solution
$$
(V_{nk}\,;~~ 0\le k\le n-1,~ n=1,\dots, N)
$$
is uniquely determined by the last row $(V_{N0},\dots, V_{N,N-1})$.
Equivalently, the corresponding Markov chain with time parameter $n$ ranging
from $N$ to $1$ is determined by the initial distribution $(\widetilde
V_{N0},\dots, \widetilde V_{N,N-1})$, which can be selected arbitrarily within
the set of all probability distributions on ${\cal E}_N$. The set ${\cal
V}^{(N)}$ of (nonnegative, normalized) solutions to such a truncated recursion
is therefore the convex hull of the arrays $V^{N\varkappa}$ ($0\le\varkappa\le
N-1$)
 which have the $N$th
row
$$
V^{N\varkappa}_{N k}= \delta_{\varkappa k}{\bigg/}\left\langle {N \atop
\varkappa} \right\rangle\,\qquad (0\le k\le N-1),
$$
where $\delta_{\varkappa k}$ is the Kronecker symbol.
The set
${\cal
V}^{(N)}$ is a $(N-1)$--dimensional simplex with extreme elements
$$
{\rm ex}({\cal V}^{(N)})= \{V^{N\varkappa}:\quad \varkappa=0,\dots,N-1\}.
$$
The  probability law $\widetilde{V}^{N\varkappa}$ corresponding to
${V}^{N\varkappa}$ rules a Markov chain which starts in state
$(N,\varkappa)\in{\cal E}_N$ at time $N$, hence the boundary of the
$N$--truncated triangle can be identified with ${\cal E}_N$.

\par
For the infinite recursion  the problem is much more complicated because there
is no obvious analogue of the ``last row'' which would provide an initial
condition for (\ref{dual}). A common recipe to obtain all solutions is the
following. Extend each $V^{N\varkappa}$ to a function on the whole set of
vertices of $\cal E$ by setting $V^{N\varkappa}_{nk}=0$ for $n>N$. Define the
{\it Martin boundary\/} \footnote{The definition corresponds to the {\it
entrance\/} boundary in \cite{Kemeny}.}
 of $\cal E$ as
$$
\partial {\cal V} := \overline{\left(\bigcup_{N, \varkappa}{\cal
V}^{N\varkappa}\right)} \setminus \left(\bigcup_{N,\varkappa}{\cal
V}^{N\varkappa}\right),
$$
where the bar means the closure in the topology of pointwise convergence of
functions on the set of vertices. Plainly, $\partial {\cal V}$ is the set of
solutions which may be obtained from truncated solutions $V^{N\varkappa}$ by
fixing some limiting regime for $\varkappa=\varkappa(N)$, as $N\to\infty$, to
secure convergence of $V^{N\varkappa}_{nk}$ for each $(n,k)$. Obviously, each
such limit is indeed a solution to (\ref{dual}), hence $\partial{\cal V}$ is a
subset of ${\cal V}$.

\par
By some well known general theory (see \cite[Ch. 1, \S1]{Kerov.book} and also
\cite{Aldous, Dynkin, Kemeny, Kerov.young}) the Martin boundary $\partial{\cal
V}$ contains the boundary $\ex(\cal V)$ (for this reason $\partial{\cal V}$ is
sometimes called the {\it maximal\/} boundary). In our proof of Theorems
\ref{thm1} and \ref{thm2} we  shall determine the Martin boundary
$\partial{\cal V}$ and then check that all its elements are actually extreme
solutions.

\par
The coincidence of boundaries is not specific for
$\cal E$, rather it holds for other number triangles and more sophisticated
graded graphs \cite{KOO}. A common reason for this phenomenon is some
law of large numbers, like the law of large numbers
for exchangeable $0\!-\!1$ random variables in the case of Pascal triangle.
On the other hand,
there are  simple examples of graded graphs for which the Martin boundary
is strictly larger than the extreme boundary \cite{Gibbs}.

\section{D--arrangements}\label{D}

Let $[n]=\{1,\dots,n\}$ and let $\Perm(n)$ be the set of permutations of $[n]$.
We write permutations $\pi\in\Perm(n)$
in the conventional one--row notation
 $\pi=\pi(1)\dots\pi(n)$
(and ignore the group structure on $\Perm(n)$).
A position $j\in[n-1]$ is said to be
a {\it descent\/} of $\pi$ if $\pi(j)>\pi(j+1)$. By $D(\pi)$ we denote the
total number of descents of $\pi$. For instance, $\pi=7356241\in\Perm(7)$ has
descents at positions $j=1,4,6$, hence $D(7356241)=3$.

\par
According to a well--known combinatorial interpretation,  the
Eulerian numbers count permutations with a given number of descents:
\begin{equation}
\Card\{\pi\in\Perm(n):\, D(\pi)=k\}=\left\langle n\atop k\right\rangle,
\end{equation}
as is easily shown by checking
that the counts satisfy the recursion \eqref{dim}
(or see \cite[Section 6.2]{GKP}).
We establish next a more delicate connection.

\par
Observe that removing $n$ from a permutation of $[n]$ yields a projection
$p_n:\Perm(n)\to\Perm(n-1)$. For instance, $3412\in\Perm(4)$ is projected to
$312\in\Perm(3)$. Clearly, the preimage of any permutation $\pi\in\Perm(n-1)$
by $p_n$ consists of exactly $n$ elements.
\begin{lemma}\label{paths}
There exists a bijection  $b_n$ between $\Perm(n)$ and standard paths in the
graph $\cal E$ of length $n$ with the following property: the path $b_n(\pi_n)$
corresponding to $\pi_n\in{\rm Perm}(n)$ passes through the vertices $(m,
D(\pi_m))\in {\cal E}_m$ $(m=1,\dots,n)$, where
$\pi_{n-1}=p_n(\pi_n),\ldots,\pi_{1}=p_{2}(\pi_{2})$ are the iterated
projections of $\pi_n$.
\end{lemma}
\begin{proof}
Choose $\pi_n\in\Perm(n)$ and let $k=D(\pi_n)$. It is readily checked
that the preimage $p_{n+1}^{-1}(\pi_n)\subset\Perm(n+1)$ consists of $k+1$
permutations with $k$ descents and of $n-k$ permutations with $k+1$ descents.
Observe that $k+1$ is the number of edges linking $(n,k)$ to
$(n+1,k)$ while $n-k$ is the number of edges linking $(n,k)$ to $(n+1,k+1)$.
It follows that if the desired bijection exists  for some
$n$ then it can be further extended to a bijection for $n+1$. The assertion
follows by induction.
\end{proof}

\vskip0.5cm
\par
The bijections $b_n$ are in no way  canonical, because we do not distinguish among
the edges linking adjacent vertices in $\cal E$. Still, the way we introduced
$b_n$'s takes care of {\it consistency} for all $n$. Indeed, let $t_n$ be the
operation of cutting off the last link in a standard path in $\cal E$ of length
$n$. Thus $t_n$  projects standard paths of length $n$ onto standard paths of
length $n-1$. The consistency of $b_n$'s amounts to the commutation relation
$t_n\circ b_n=b_{n-1}\circ p_n$, which holds  for all $n\geq 2$.

\par
Let $\A=\varprojlim\Perm(n)$ be the inverse limit \footnote{Another inverse
limit, the  space of virtual permutations, appears in \cite{KOV}.} of the
finite permutation spaces $\Perm(n)$ with respect to $p_n$'s. Elements of $\A$
are infinite sequences $(\pi_n)$ of {\it consistent} permutations
 $\pi_n\in\Perm(n)$ ($n=1,2,\ldots$), meaning that, for each $n\ge2$,
 $p_n(\pi_n)=\pi_{n-1}$. In extension of Lemma \ref{paths} we have the following
 corollary.

\begin{corollary}\label{paths2}
The consistent sequence of bijections $(b_n)$ defines a bijection
between
$\A$ and the set of infinite standard paths in
 $\cal E$. The bijection has the property that
the path corresponding to
$(\pi_n)\in\A$  passes
through the vertices $(n, D(\pi_n))$,
$n=1,2,\ldots$.
\end{corollary}

\vskip0.5cm
\par
With each permutation $\pi\in\Perm(n)$ we associate a total order $\lhd$ on the
set $[n]$, in which $\pi(1)\lhd\dots\lhd\pi(n)$. Likewise, every element
$(\pi_n)\in\A$ determines a total order on the set ${\mathbb N}=\{1,2,\dots\}$
such that, for each $n$, the total  order restricted to the subset
$[n]\subset{\mathbb N}$ is the one given by $\pi_n$. Conversely, any total
order on ${\mathbb N}$ can be obtained in this way, from some element of $\A$.
For this reason, we call the elements of $\A$ {\it arrangements\/} and identify
them with the total orders on $\Nat$.Two obvious examples of arrangements are
 the standard order $1\lhd2\lhd3\lhd\dots$ and  the inverse order
$\dots\lhd3\lhd2\lhd1$; the corresponding paths in $\cal E$ go along the left
side of the Euler triangle and along its right side, respectively.

\par
As a projective limit of finite sets, $\A$ is a compact topological space.
Given a probability measure $P$ on $\A$ we can speak of a {\it random}
arrangement $\Pi=(\Pi_n)$, where $\Pi_n\in \Perm(n)$ are consistent random
permutations, such that the law $P_n$ of $\Pi_n$ is the pushforward of $P$ by
the canonical projection $\A\to\Perm(n)$. Conversely, by Kolmogorov's measure
extension theorem, each sequence of distributions $(P_n)$ determines  a unique
random arrangement, provided the sequence is consistent with respect to all
projections $p_n$.

\par
The random arrangements relevant to our discussion have one special property
of sufficiency.

\begin{definition} We say that a random arrangement $\Pi=(\Pi_n)$ is a {\it
D--arran\-ge\-ment\/} if for every $n=1,2\dots$ and  $\pi_n\in\Perm(n)$ the
probability of the event $\Pi_n=\pi_n$ depends on the couple $(n,D(\pi_n))$
only.
\end{definition}
\noindent
That is to say, for a D--arrangement $\Pi=(\Pi_n)$, the number of
descents is a {\it sufficient statistic\/}: conditionally given $D(\Pi_n)=k$
the distribution of $\Pi_n$ is uniform on the set of  permutations of $[n]$
with $k$ descents, for each $n$ and $k$.

\par
Two trivial examples of D--arrangements are the nonrandom arrangements given by
the standard order and the inverse order. The corresponding measures on $\A$
are the Dirac masses at  points $(\pi_n)=(1\dots n)\in\A$ and
$(\pi_n)=(n\dots1)\in\A$, respectively. Notice that these two are the only
Dirac measures on $\A$ corresponding to D--arrangements.

\par
More substantial example is the random arrangement $\Pi$
for which every $\Pi_n$ has  uniform distribution on $\Perm(n)$. This is the
only {\it exchangeable\/} random arrangement, whose probability law is
invariant under arbitrary permutations of the set ${\mathbb N}$.

\par
Now, Corollary \ref{paths2} implies:

\begin{lemma}\label{bijection}
The formula
$$
P_n(\pi_n)=V_{n,D(\pi_n)}, \qquad n=1,2,\dots, \quad \pi_n\in\Perm(n)
$$
defines an affine isomorphism $V=(V_{nk})\leftrightarrow P=(P_n)$
between
$\cal V$ and the set of probability laws for D--arrangements.
\end{lemma}
\vskip0.5cm

\noindent Equivalently, in terms of quantities $\widetilde V_{nk}$ and random
paths in $\cal E$ corresponding to D--arrangements, $\widetilde V_{nk}$ is the
probability that a random infinite path (with distribution $P$) will pass
through the vertex $(n,k)$. In the sequel we will not distinguish between
solutions to (\ref{dual}) and random D--arrangements.

\section{Bucket sorting}\label{B}

Here we use the correspondence of Lemma \ref{bijection} for constructing a
family of solutions $V\in \cal V$. The following algorithm, called  {\it bucket
sorting\/}, exploits a multinomial distribution and is a simplest of the
algorithms of this kind, widely known in computer science \cite{Mahmoud} and
dynamical systems \cite{BD, Lalley}.

\par
Fix $\varkappa\in \{0,1,2,\dots\}$ and imagine $\varkappa+1$ buckets
arranged in some order. Suppose each of the numbers $1,2,\dots\in{\mathbb N}$ is sent
to one of the buckets with equal probabilities $(\varkappa+1)^{-1},
\dots,(\varkappa+1)^{-1}$, independently of the other numbers.
 For each $n$ this yields a random allocation
of integers $1,\dots,n$ in the buckets.
Arranging the integers within each
bucket in increasing order and
putting the resulting sequences together
(in the order of the buckets), the allocation of $n$ integers is transformed
into a random permutation $\Pi^\varkappa_n$ of
$[n]$. By the construction,  $\Pi^\varkappa_n$ has at most $\varkappa$ descents.

\begin{lemma}\label{bucket}
The infinite sequence $\Pi^\varkappa=(\Pi^\varkappa_n)$ produced by the bucket
sorting  is a D--arrangement. The corresponding array
$W^\varkappa=(W^\varkappa_{nk})\in\cal V$ is given by formula
\begin{equation}\label{bucket1}
W^\varkappa_{nk}={n+\varkappa-k\choose n}{\bigg/}(\varkappa+1)^n.
\end{equation}
\end{lemma}

\begin{proof}
By the very construction, the random permutations $\Pi^\varkappa_n$ are
consistent with respect to the projections $p_n$, hence $\Pi^\varkappa$ is a
random arrangement. Given $\pi_n\in\Perm(n)$, let us compute the probability of
the event $\Pi^\varkappa_n=\pi_n$. The total number of possible allocations of
$1,\dots,n$ into buckets equals $(\varkappa+1)^n$, and all of them are equally
likely. Thus, it suffices to compute the number of the allocations resulting in
$\pi_n$. Any such allocation is determined by a partition of the sequence
$\pi_n=\pi_n(1)\dots\pi_n(n)$ into $\varkappa+1$ consecutive fragments (some of
which can be empty), and any such partition can be encoded by placing
$\varkappa$ vertical bars separating the fragments. Observe that for each
descent $j\in[n-1]$ of the permutation $\pi_n$ at least one bar has to be
placed between $\pi_n(j)$ and $\pi_n(j+1)$. For $k=D(\pi_n)$ we see that $k$
positions of bars are fixed by the descents, so that the allocation is actually
determined by the remaining $\varkappa-k$ bars, which can be placed
arbitrarily. Since the bars are indistinguishable, the number of possibilities
equals ${n+\varkappa-k\choose n}$. Thus, the probability of $\pi_n$ is given by
the right--hand side of \eqref{bucket1}. Since this expression depends only on
$k=D(\pi_n)$, we conclude that $\Pi^\varkappa$ is a D--arrangement and
\eqref{bucket1} is the corresponding element of $\cal V$.
\end{proof}

\begin{remark}
The fact that  formula \eqref{bucket1} determines a solution to \eqref{dual}
amounts to a binomial identity which is easy to check directly:
$$
(\varkappa+1)\,{n+\varkappa-k\choose n}= (k+1)\,{n+\varkappa+1-k\choose n+1}+
(n-k)\,{n+\varkappa-k\choose n+1}\,,
$$
while the total probability rule  \eqref{normalisation} becomes
$$
(\varkappa+1)^n=\sum_{k=0}^{n-1}\left\langle{n\atop k}\right\rangle
{n+\varkappa-k\choose n},
$$
which is equivalent to Worpitzky's identity \cite[(6.37)]{GKP}.
\end{remark}

Recall that in Section \ref{M} we introduced arrays
$V^{N\varkappa}=(V^{N\varkappa}_{nk})\in{\cal V}^{(N)}$ solving the
``$N$--truncated'' version of recursion \eqref{dual}.

\begin{lemma}\label{approx}
Fix $\varkappa\in \{0,1,\dots\}$, let $\Pi^\varkappa=(\Pi^\varkappa_n)$ be the
D--arrangement resulting from the bucket sorting, and let
$W^\varkappa=(W^\varkappa_{nk})\in{\cal V}$ stand for the corresponding array.
Then $V^{N\varkappa}$ converge to $W^\varkappa$, that is
$$
\lim_{N\to\infty}V^{N\varkappa}_{nk}=W^\varkappa_{nk}\,~~~(n=1,2,\ldots;\,
0\leq k\leq n-1).
$$
\end{lemma}

\begin{proof}
It is more convenient to deal with quantities
$$
\widetilde V^{N\varkappa}_{nk}=\left\langle{n\atop k}\right\rangle
V^{N\varkappa}_{nk}\,, \qquad \widetilde W^\varkappa_{nk}=\left\langle{n\atop
k}\right\rangle W^\varkappa_{nk}\,.
$$
For fixed $n$, the vectors
$$
(\widetilde V^{N\varkappa}_{n0}, \dots, \widetilde V^{N\varkappa}_{n, n-1})
\quad \text{and} \quad (\widetilde W^\varkappa_{n0}, \dots, \widetilde
W^\varkappa_{n,n-1})
$$
are the distributions at time $n$ of the $N$--step Markov chain (introduced in
Section \ref{M}) whose initial distribution at time $N$ is
\begin{equation}\label{initial}
(\widetilde V^{N\varkappa}_{N0}, \dots, \widetilde V^{N\varkappa}_{N, N-1})
\quad \text{or} \quad (\widetilde W^\varkappa_{N0}, \dots, \widetilde
W^\varkappa_{N,N-1}),
\end{equation}
respectively.

\par
Recall that both $\widetilde V^{N\varkappa}_{Nk}$ and $\widetilde
W^\varkappa_{Nk}$ vanish for $k>\varkappa$ and, moreover,
\begin{equation}\label{deltav}
\widetilde
V^{N\varkappa}_{Nk}=\delta_{\varkappa k}.
\end{equation}
 We claim that it suffices to prove
the limit relation
\begin{equation}\label{limit}
\lim_{N\to\infty}\widetilde W^\varkappa_{N\varkappa}=1.
\end{equation}
Indeed, since (\ref{initial}) are probability distributions,
(\ref{deltav}) and (\ref{limit}) imply that the total variance distance between them
goes to $0$ as $N\to\infty$,
which implies the assertion of the lemma.

\par
To prove \eqref{limit} we turn to $\Pi^\varkappa$ and observe that $\widetilde
W^\varkappa_{N\varkappa}$ is just the probability for the random permutation
$\Pi^\varkappa_N$ to have the maximal possible number of descents $\varkappa$.
In terms of the random allocation of the numbers $1,\dots,N$, this means that
all buckets are nonempty and the largest number in each bucket (except the last
bucket) is larger than the smallest number in the next bucket. If this were not
the case, all the numbers in one of the buckets were smaller than those in the
next bucket. Elementary estimates which we postpone to the proof of Lemma
\ref{char} show that the probability of such an event tends to $0$ as
$N\to\infty$, which yields  \eqref{limit}.
\end{proof}

\begin{remark}
By virtue of the explicit formula \eqref{bucket1}, the limit relation
\eqref{limit} is equivalent to an asymptotic relation for the Eulerian numbers:
$$
\lim_{N\to\infty} \left\langle N \atop \varkappa \right\rangle \sim
(\varkappa+1)^N \qquad \text{for fixed $\varkappa=0,1,\dots$}.
$$
This relation can be readily checked directly. For instance, it follows from
formula \eqref{explformula} below.
\end{remark}

\par
Lemma \ref{approx} shows that the family $\{W^\varkappa\}$ is contained in the
Martin boundary. Actually, a stronger claim holds: all $W^\varkappa$'s are
extreme. We show this in Lemma \ref{extreme} below. But first we will prove
the law of large numbers for $\Pi^\varkappa$.

\begin{lemma}\label{char}
Fix $\varkappa=\{0,1,\dots\}$. We have
$$
\lim_{n\to\infty}D(\Pi^\varkappa_n)=\varkappa \qquad\text{\rm with probability
1}.
$$
Moreover, this property of $\Pi^\varkappa$ is characteristic.
\end{lemma}
\begin{proof}
Suppose there are just two buckets, $\varkappa=2$. Then $D(\Pi^2_{n})<2$ means
that for some $m\in [n]$ the integers $1,\ldots,m$ fall in the first bucket,
and $m+1,\ldots,n$ in the second, which is an event of probability $(n+1)/2^n$.
Since the series of these probabilities converges, the Borel--Cantelli lemma
yields the claim. The general case $\varkappa>2$ is reduced to the estimate in
the case $\varkappa=2$ by focussing on two consecutive buckets and using
elementary  large deviation bounds for Bernoulli trials to show that the chance
for  less than, say, $n/\varkappa$ integers in both buckets goes to $0$
exponentially fast with $n$. The uniqueness follows as in Lemma \ref{approx}.
\end{proof}

\begin{lemma}\label{extreme}
Elements $W^\varkappa\in\cal V$ resulting from the bucket sorting are
extreme.
\end{lemma}

\begin{proof}
If $W^\varkappa$ is a mixture of some $V^1,V^2\in {\cal V}$ then
by the first assertion of Lemma \ref{char} the arrangements corresponding to
$V_1$ and $V_2$ must satisfy the same law of large numbers as $W^\varkappa$.
But then by the second assertion of the lemma $V_1=V_2=W^\varkappa$.
Hence $W^\varkappa$ is extreme.
\end{proof}

\section{Proofs of the main results}
\label{P}

We start with reducing the set of parameters needed to determine a generic
solution $V\in {\cal V}$.
\begin{lemma}\label{leftside}
The sequence $(V_{n0})$
uniquely determines
 $V\in\cal V$.
\end{lemma}

\begin{proof}
Writing \eqref{dual} as
$$
V_{n+1,k+1}=\frac1{n-k}\,V_{nk}-\frac{k+1}{n-k}\,V_{n+1,k}\, \qquad (0\le k\le
n-1),
$$
we see that for each $k\in\{0,1,\dots\}$  the sequence
$(V_{n,k+1}\,:\, n=k+2,k+3,\dots)$ is uniquely determined by the sequence
$(V_{nk}\,:\, n=k+1,k+2,\dots)$. Induction in $k$ ends the proof.
\end{proof}
\vskip0.5cm

\par
Consider the truncated arrays $V^{N\varkappa}=(V^{N\varkappa}_{nk})\in{\cal
V}^{(N)}$ with parameters $N$ and $\varkappa$, $0\le\varkappa\le N-1$,
introduced in Section \ref{B}. As in the proof of Lemma \ref{approx} we
introduce the modified array $\widetilde V^{N\varkappa}=(\widetilde
V^{N\varkappa}_{nk})$ with $\widetilde V^{N\varkappa}_{nk}=\left\langle n\atop
k\right\rangle V^{N\varkappa}_{nk}$, and we recall that the $n$th row of
$\widetilde V^{N\varkappa}$ is the distribution at time $n$ of the  Markov
chain started at time $N$ from the vertex $(N,\varkappa)$. Since
$V^{N\varkappa}_{n0}=\widetilde V^{N\varkappa}_{n0}$, the quantity
$V^{N\varkappa}_{n0}$ equals the probability of the event that the Markov chain
will pass through vertex $(n,0)$.

\par
Furthermore, there is a monotonicity property analogous to that
of generalized Stirling triangles in \cite{Gibbs}.
\begin{lemma}\label{monotone}
For  fixed $n\in \{1,\dots,N\}$, the coordinate $V^{N\varkappa}_{n0}$ does not
increase as $\varkappa$ varies from $0$ to $N-1$.
\end{lemma}

\begin{proof}
We employ a simple coupling argument. Given two numbers $\varkappa<\varkappa'$,
consider two Markov chains which start from vertices $(N,\varkappa)$ and
$(N,\varkappa')$, respectively. We settle both chains on a common probability
space assuming that the jumps are independent as long as the  trajectory of the
first chain does not intersect the trajectory of the second chain, but once the
trajectories meet, they merge. The merge  does not affect the marginal law of
each of the chains, since both are directed by the same transition
probabilities. The key property of the coupling is that each trajectory of the
first chain remains on the left of the trajectory of the second chain, before
the trajectories merge. Observe now that after reaching the left side of the
Euler triangle, a trajectory can only process along this side. Consequently, if
a trajectory of the second chain passes through $(n,0)$ then the trajectory of
the first chain reaches the left side of the triangle at some time $n'\ge n$,
hence passes through $(n,0)$, too. Thus, the chance for the first chain to pass
through $(n,0)$ is not less than that for the second chain. This proves the
desired inequality $V^{N\varkappa}_{n0}\ge V^{N\varkappa'}_{n0}$.
\end{proof}
\vskip0.5cm

\par
We proceed with the proof of Theorems \ref{thm1} and \ref{thm2}.
Our strategy is to determine first the Martin boundary
by directly identifying all solutions $W$ that
appear as limits  of
arbitrary sequences of the form $V^{N,\varkappa(N)}$.

\par
Assume first that $\varkappa(N)=\varkappa$ with some fixed $\varkappa$, for all
$N$ large enough. Then, by Lemma \ref{approx}, the sequence
$V^{N,\varkappa(N)}$ converges to the array $W^\varkappa$ given by formula
\eqref{bucket1} which we also display here for reader's convenience:
\begin{equation}\label{1}
W^\varkappa_{nk}={n+\varkappa-k\choose n}{\bigg/}(\varkappa+1)^n.
\end{equation}

\par
Next, assume that $\varkappa(N)=N-1-\varkappa$, where $\varkappa$ is
fixed. Observe that this limit regime is reduced to the preceding one by
application of the symmetry $(n,k)\to(n,n-1-k)$ of the Euler triangle $\cal E$.
Therefore, in this case the sequence converges to the array $\widehat W^\varkappa$
with components
\begin{equation}\label{2}
\widehat W^\varkappa_{nk}={\varkappa+k+1\choose n}{\bigg/}(\varkappa+1)^n.
\end{equation}
The random D-arrangement corresponding to $\widehat W^\varkappa$ can be produced by
the obvious analogue of bucket sorting in which integers within each bucket are arranged
in {\it decreasing} order.

\par
Further on, from \eqref{1} and \eqref{2} it is readily seen that there exists
the limit
$$
\lim_{\varkappa\to\infty}W^{\varkappa}=\lim_{\varkappa\to\infty}\widehat
W^{\varkappa}=W^\infty
$$
with components
$$
W^\infty_{nk}=\frac1{n!}.
$$
Clearly, $W^\infty\in\cal V$.

\par
Now we claim that if both $\varkappa(N)$ and $N-1-\varkappa(N)$ go to infinity
then the sequence $V^{N,\varkappa(N)}$ converges to $W^\infty$. To that end,
observe that a general bound
\begin{equation}\label{bound}
0\le V_{nk}\le \left\langle n\atop k\right\rangle^{-1},
\end{equation}
which follows from Lemma \ref{norm}, holds for all $V\in\cal V$ and implies
that $\cal V$ is compact in the product topology. By the compactness, passing
if necessary to a subsequence of $(V^{N,\varkappa(N)})$ we can always achieve
convergence to some $W$, hence it is enough to show that $W=W^\infty$, and by
Lemma \ref{leftside}, this is further reduced to showing that
$W_{n0}=W^\infty_{n0}=\frac1{n!}$. For any fixed $\varkappa$ we have
$$
\varkappa\le \varkappa(N)\le N-1-\varkappa
$$
for large $N$. Applying Lemma \ref{monotone} we obtain the bound
$$
W^\varkappa_{n0}\ge W_{n0}\ge \widehat W^\varkappa_{n0}\, \qquad (n=1,2,\dots).
$$
Now, sending $\varkappa$ to infinity we conclude that $W_{n0}=\frac1{n!}$, as wanted.

\par
We have shown that the Martin boundary consists of the elements
$$
W^\varkappa \quad (\varkappa=0,1,2,\dots), \quad \widehat W^\varkappa \quad
(\varkappa=0,1,2,\dots), \quad {\rm and}\quad  W^\infty.
$$
Comparison with formulas of Theorem \ref{thm1} shows that these are exactly the
arrays $W(\theta)$ with $\theta>\frac12$, $\theta<\frac12$, and
$\theta=\frac12$, respectively. A remarkable fact emerges:  the single entry $(2,0)$
distinguishes all these arrays.

\par
By Lemma \ref{extreme}, the arrays $W^\varkappa$ are extreme. By symmetry, the
arrays $\widehat W^\varkappa$ are extreme, too. To finish the proof of the
theorems it remains to check that $W^\infty=W(\frac12)$ is extreme. We postpone this
to the next section.

\par
The Corollary \ref{uniqueness} follows from general results.
Since the space ${\cal V}$ is  compact, metrizable and separable,
 the well--known Choquet theorem
\cite[\S3]{Phelps} implies that each solution $V\in\cal V$ may be represented as a convex
mixture of the extreme solutions $W\in {\rm ex}({\cal V})$.
A simple general argument shows that
$\cal V$ is a Choquet simplex (that is, the cone generated by $\cal V$ is a
lattice cone), see, e.g., \cite[Lemma 9.3]{Olshanski2}.
The uniqueness of representation now follows from another
Choquet's theorem, see \cite[\S9]{Phelps}.

\section{End of proof: extremality of the exchangeable arrangement}
\label{EOP}
In this section  $P$ denotes the probability measure on $\A$ corresponding to the array
$W^\infty=W(\frac12)\in\cal V$. The characteristic property of $P$ is that, for
each $n$, the image $P_n$ under the natural projection $\A\to\Perm(n)$ is the
uniform measure assigning to all permutations $\pi_n\in\Perm(n)$ equal weights
$\frac1{n!}$. Let $\Pi$ be the random arrangement with law $P$.
This is the  exchangeable random arrangement,
invariant under the
natural action on the space $\A$ of
permutations of
$\mathbb N$.

\par
The following useful construction of $\Pi$ is found in \cite{Aldous}.
Let $X_1,X_2,\dots$ be  independent random variables,  with  uniform
distribution on the unit interval $[0,1]$. The  $X_j$'s are pairwise  distinct
with probability one. Define a random total order on $\mathbb N$ by the rule
$i\lhd j$ if $X_i< X_j$.  Clearly, for each $n$, the resulting random
permutation $\Pi_n$ of $[n]$ depends only of $X_1,\dots,X_n$ and, by
exchangeability of $X_j$'s, $\Pi_n$ is uniformly distributed on $\Perm(n)$.

\par
Known moments of $D(\Pi_n)$ follow easily from this realization.
\begin{lemma}\label{LL}
Let $\Pi_n$ be the uniform random permutation of $[n]$.
The random variable $D(\Pi_n)$ has mean $(n-1)/2$ and
variance $(n-1)/12$.
\end{lemma}

\begin{proof}
Clearly, $D(\Pi_n)$ equals the number of descents in the random sequence
$X_1,\dots,X_n$, that is, the number of indices $i\in[n-1]$ such that
$X_i>X_{i+1}$. Therefore, denoting by $\chi_i$ the indicator of the event
$X_i>X_{i+1}$ we have
$$
D(\Pi_n)=\sum_{i=1}^{n-1}\chi_i\,.
$$
The result easily follows from the relations
$$
{\mathbb E}(\chi_i)={\mathbb E}(\chi_i^2)=\tfrac12, \qquad {\mathbb
E}(\chi_i\chi_{i+1})=\tfrac16\,, \qquad {\mathbb E}(\chi_i\chi_j)=0
\quad(|i-j|\ge2).
$$
\end{proof}

Since both the mean and the variance exhibit a linear  growth, standard
application of Chebyshev's inequality gives:

\begin{corollary}\label{cor}
Under the uniform distribution,
$D(\Pi_n)/(n-1)\to\frac12$ in probability.
\end{corollary}

We have now all tools to finish the argument of Section \ref{P} by showing that
$W^\infty$ is extreme. Assume the contrary, then the boundary $\ex(\cal V)$
reduces to $\{W^\varkappa\}\cup\{\widehat W^\varkappa\}$ and hence $P$ can be
written as a convex combination of the measures $P^\varkappa$ and $\widehat
P^\varkappa$ (the laws of $W^\varkappa$ and $\widehat W^\varkappa$),
$\varkappa=0,1,2,\dots$. By Corollary \ref{cor}, there exists a sequence of
numbers $n_1<n_2<\dots$ such that $D(\Pi_{n_i})/(n-1)\to\frac12$ almost surely.
On the other hand the same ratio goes to  $0$ or $1$ under the distribution
$P^\varkappa$ or $\widehat P^\varkappa$, respectively. This leads to a
contradiction, so the proof is complete.

\section{Concluding remarks}

\subsection{Permutations with descent--set statistic}

We were led to consider the Eulerian triangle $\cal E$ in connection with a
larger graded graph $\cal Z$ of zigzag diagrams \cite{GO}. With edge
multiplicities taken into account, both graphs have the same path spaces, but
$\cal Z$ has more vertices and much more rich branching. The boundary problem
for $\cal Z$ amounts to describing all random arrangements $\Pi=(\Pi_n)$ with
the property that the distribution of each $\Pi_n$ is uniform conditionally
given the set of descent positions in $\Pi_n$. In \cite{GO} we established that
the distribution of a random total order determined by such $\Pi$ must be
spreadable, that is invariant under increasing mappings ${\mathbb N}\to{\mathbb
N}$. D--arrangements are the simplest of this kind, and the extreme
D--arrangements we described here are also extreme solutions to the boundary
problem for $\cal Z$.

\par
Analogous connection exists between Kingman's graph $\cal K$ of partitions and
the Stirling triangle $\cal S$ of the first kind \cite{Gibbs}. The relevant
random objects are exchangeable partitions of ${\mathbb N}$ and a smaller class
of partitions which have the number of blocks as sufficient statistic. In that
case the situation is more interesting than the one for $\cal Z$ and $\cal E$:
extremes solutions to the boundary problem for $\cal S$ (the celebrated Ewens
partition structures) are decomposable along the boundary of $\cal K$, with the
mixing measure being the remarkable
 Poisson--Dirichlet distribution \cite{CSP}.

\subsection{A problem of moments}

Corollary \ref{uniqueness} and (\ref{GrishinaF}) tell us that every solution
$V\in {\cal V}$ satisfies
$$
V_{n0}=\sum_{\theta\in \Theta} p(\theta)
\,\,\prod_{i=0}^{n-1}\,{1+(2\theta-1)i\over 1+i}
$$
for some unique probability distribution $p$ on the parameter set $\Theta$. An
inverse problem asks one to characterize all sequences $(V_{n0})$ with
$V_{10}=1$ which can be represented in this form. An answer is suggested by the
argument in Lemma \ref{leftside} which says that  there is a linear operator
$\nabla: (V_{n0})\mapsto (V_{nk})$ which maps an arbitrary sequence to a
solution of (\ref{dual}). So the necessary and sufficient condition for
representability is that $\nabla$ applied to $(V_{n0})$ produces a nonnegative
array.

\par
The analogous question for Pascal's triangle is  the Hausdorff moment problem,
with kernel $\theta^n$, where $\theta$ ranges in $[0,1]$. In this classical
case the analogue of $\nabla$ associates with each sequence the array of its
iterated differences, whose positivity is Hausdorff's condition called total
monotonicity.

\subsection{Remarks on the uniform case}
The Eulerian numbers are given by the formula
\begin{equation}\label{explformula}
\left\langle {n \atop k}
\right\rangle=\sum_{j=0}^{k} (-1)^j {n+1\choose j}(k+1-j)^n\,,
\end{equation}
which compared with Laplace's formula in \cite[Section 1.9]{Feller} shows that
$${\rm Prob}(k\leq Y_1+\ldots+Y_n<k+1)=\left\langle{n\atop k}\right\rangle/n!$$
for $Y_1,Y_2,\dots$ independent  random variables with
uniform distribution on $[0,1]$.
The following explanation
 of this coincidence is borrowed from
\cite[p. 296]{Pitman}. For $x>0$ let $x=\lfloor x\rfloor+\{x\}$ be the
decomposition of $x$ into integer and fractional parts, with
 $0\le\{x\}<1$.
Consider
$$
S_j=Y_1+\dots+Y_j, \qquad X_j=\{S_j\},
$$
then $X_1, X_2,\ldots$ are also independent, uniform on $[0,1]$.
Observe that $X_{j}>X_{j+1}$ each time $\lfloor S_j \rfloor <\lfloor S_{j+1}\rfloor
 =\lfloor S_j\rfloor+1$ and recall the discussion preceding Lemma \ref{LL}.

\par
Improving upon Corollary \ref{cor},
we see that the convergence $D(\Pi_n)/(n-1)\to\frac12$  holds with probability
1. The central limit theorem applied to $S_n$'s entails that the distribution
of $D(\Pi_n)$ is asymptotically Gaussian. This connection with sums of random
variables has been a starting point for many fine  results on  descents in
uniform permutation. See \cite{Oshanin} for recent development and references.

\def\cprime{$'$} \def\polhk#1{\setbox0=\hbox{#1}{\ooalign{\hidewidth
\lower1.5ex\hbox{`}\hidewidth\crcr\unhbox0}}} \def\cprime{$'$}
\def\cprime{$'$} \def\cprime{$'$}
\def\polhk#1{\setbox0=\hbox{#1}{\ooalign{\hidewidth
\lower1.5ex\hbox{`}\hidewidth\crcr\unhbox0}}} \def\cprime{$'$}
\def\cprime{$'$} \def\polhk#1{\setbox0=\hbox{#1}{\ooalign{\hidewidth
\lower1.5ex\hbox{`}\hidewidth\crcr\unhbox0}}} \def\cprime{$'$}
\def\cprime{$'$} \def\cydot{\leavevmode\raise.4ex\hbox{.}} \def\cprime{$'$}
\def\cprime{$'$} \def\cprime{$'$} \def\cprime{$'$}

\end{document}